\documentclass{article}
\newcommand{\proof}{\noindent {\bf Proof: }}
\newcommand{\note}{\noindent {\bf Notation: }}
\newcommand{\remark}{\noindent {\bf Remark:}}
\newcommand{\corollary}{\noindent {\bf Corollary:}}
\newcommand{\ws}{\hspace{4pt}}
\newtheorem{theorem}{Theorem}

\newtheorem{lemma}{Lemma}
\newtheorem{defi}{Definition}
\usepackage[]{amssymb}

\begin{document}
\title{Biorthonormal Systems in Freud-type Weighted Spaces with Infinitely Many Zeros - An Interpolation Problem}
\author{\'Agota P. Horv\'ath \footnote{supported by Hungarian National Foundation for Scientific Research, Grant No. T049301.\newline Key words:Complete and minimal system, infinite interpolation, infinite linear equation system, Freud weight. \newline  2000 MS Classification: 42A65, 41A05}}\date{}
\maketitle
\begin{center}\small{Department of Analysis, Budapest University of Technology and Economics}\end{center}
\begin{center}\small{H-1521 Budapest, Hungary  }\end{center}
 \begin{center}\small{e-mail: ahorvath@renyi.hu}\end{center}

\begin{abstract}In a Freud-type weighted ($w$) space, introducing another weight ($v$) with infinitely many roots, we give a complete and minimal system with respect to $vw$, by deleting infinitely many elements from the original orthonormal system with respect to $w$. The construction of the conjugate system implies an interpolation problem at infinitely many nodes. Besides the existence, we give some convergence properties of the solution.
\end{abstract}

\section{Introduction}

In several related problems occurs the claim of constructing biorthonormal systems in certain Banach spaces. The initial investigations of eg. R. P. Boas and H. Pollard, A. A. Talalyan, M. Rosenblum, and B. Muckenhoupt resulted the development of eg. $A_p$-weights, the theory of multiplicative completion of set of functions, estimations on certain norms of Poisson integrals (\cite{boa}, \cite{ta}, \cite{ro}, \cite{mu}). Furter results were given eg. on completion (\cite{pz}, \cite{kz}, \cite{k2}), solving Dirichlet's problem with respect to boundary functions with singularities (\cite{k}, \cite{kh}), constructing A-basis (basis for Abel-summability) in some Banach spaces (\cite{gha}). 

The reason, why we are interested in constructing complete and minimal systems, (We say that $\{\varphi_n\}$ is minimal in a Banach space $B$, if there is a conjugate system in the dual space $\{\varphi_n^{*}\}\subset B^{*}$, such that $\varphi_n^{*}(\varphi_m) = \delta_{n,m},$ and $\{\varphi_n\}$ is complete, when for a $g \in B^{*}, g(\varphi_n)=0, n\in \mathbb{N}$ implies that $g=0$.) is the following theory of S. Banach (\cite{ban}):

\newpage

\noindent {\bf Theorem X}

{\it A system $\{\varphi_n\}_{n=n_0}^{\infty}$ is A-basis in the space $L^p_{vw}$ $(1< p < \infty)$ (with some weight function $vw$) if and only if it is a complete and minimal system in $L^p_{vw}$, and there is a constant $c = c(p)$ such that
$$\sup_{0\leq r < 1}\left\|\sum_{n=n_0}^{\infty}r^na_n(f)\varphi_n\right\|_{vw,p} \leq c\|f\|_{vw,p}$$
where $a_n(f)= \varphi_n^{*}f = \int_{{\bf R}}f\varphi_n^{*}v^2w^2$.}

So according to Banach's theorem, if a complete and minimal system is given in a Banach space, then for proving that this system is an $A$-basis, it is enough to show that the norm of the Poisson integral is bounded by the norm of the function.

On the language of weighted spaces on the real line, the common idea of the above-mentioned investigations is the following: there is a complete orthonormal system $\{e_n\}$ with respect to a weight $w>0$ (sometimes $w\equiv$ 1) on a finite or infinite interval $I$, and $v$ is another weight on $I$ with some zeros. Removing some elements of $\{e_n\}$ (which omission depends on the roots of $v$), a complete and minimal system can be constructed in a weighted space with respect to $vw$, wich means in other terminology, that the residual system can be multiplicatively trasformed into a basis. If the number of the roots ($M$) of $v$ is finite (with multiplicity), then the biorthonormal system $\{\varphi_n,\varphi_m^{*}\}$ will be the following: 
$$\{\varphi_n\}= \{e_n\}\setminus \{e_{k_1},\dots, e_{k_M}\},$$ 
and the elements of the conjugate system will be 
$$\varphi_m^{*}= \frac{\varphi_m-\sum_{i=1}^Ma_{im}e_{k_i}}{v^2}.$$
Here the denominator has some zeros, so roughly speaking, the nominator has to be zero at the same points with the same multiplicity, wich results an interpolation problem. Eg., if $\{e_n\}=\{p_n\}$ is the orthogonal polynomial system on an interval, with respect to a weight, then the linear combination of the first $M$ elements, which is a polynomial with degree $M-1$, interpolates the residual elements at the zeros of $v$ (\cite{gha}, \cite{kh}, \cite{k}). Generally, in the finite case we get a finite linear equation system, and if it has a unique solution, the biorthonormal system is complete and minimal. 

The question is the following: what can we do, if $v$ has infinitely many zeros? Following the same chain of ideas, we have to remove infinitely many elements of the original basis such that at the roots of $v$ the elements of the residual system can be interpolated by an infinite combination of the removed elements. (Naturally, we can not omit the first infinitely many elements of the original system.) We have to solve an infinite interpolation problem, which implies an infinite linear equation system. That is besides the solvability of the equation system and the unicity of the solution, the convergence of the solution (in some sense) is also a problem.

We will carry out this type of investigations on the real line, when the "outer" weight will be a Freud weight. The ideal situation would be that for an almost arbitrary root system (eg. when it has no finite accumulation point) of an "inner" weight $v$, which does not grows too quick at infinity, one could give a good omission system, but at present we are unable to state any result in this respect.

 Supposing some polynomially uniform growing property of the choosing function, we will be able to construct a point system, which will be the zeros of $v$, and an omission system step by step, with which the residual system will be complete and minimal.
 Furthermore we can apply a finite section method (\cite{gro}) to get the numerical solution of the infinite equation system.

\section{Definitions, Notations, Result}

At first we define Freud weights as generally as we will use in this paper.

\begin{defi} \cite{lelu} $w(x) = e^{-Q(x)}$ is a Freud weight, if $Q:\mathbb{R} \longrightarrow\mathbb{R}$ is even, continuous in $\mathbb{R}$, $Q(0) =0$, $Q^{''}$ is continuous in $[0, \infty)$, and $Q^{'} > 0 $ in $(0, \infty)$. Furthermore, assume that for some $A, B > 1,$
\begin{equation} A \leq \frac{(d/dx)(xQ^{'}(x))}{Q^{'}(x)} \leq B, \ws\ws  x \in (0,\infty) \end{equation} \end{defi}

\noindent {\bf Notation:}

\noindent (1) For a Freud weight $w$ we will denote by $p_n(w)=p_n$ the $n^{th}$ orthonormal polynomial on the real line, with respect to $w^2$.

\noindent (2) $w$ is a weight function 
\begin{equation} f \in L_w^p \ws\ws  \mbox{iff} \ws\ws fw \in L^p.\end{equation}
If $f \in L_w^p$ and $g \in L_w^q$, where $\frac{1}{p} + \frac{1}{q} = 1$ then let us denote by: 
\begin{equation} <f,g> = \int_{\mathbb{R}} fgw^2 \end{equation}

\vspace{2mm}

 After the definition of the external weight we give the form of that part of the weight function which is responsible for the inner roots.
The definition below is based on the Lemma 1.1 of J. Szabados \cite{joska}

\begin{defi} Let $X := \{ x_1, x_2, \dots\} \subset \mathbb{R}$, $0<|x_1| \leq |x_2|\leq\dots $ be a point system on the real line, and let $M:= \{m_1, m_2, \dots \}\subset \mathbb{R}_+$ be a collection of positive numbers. If there exists a nonnegative number $\varrho \geq 0$ such that 
\begin{equation} \sum_{j=1}^{\infty}\frac{m_j}{|x_j|^{\varrho+\varepsilon}} < \infty, \mbox{ \ws but \ws } \sum_{j=1}^{\infty}\frac{m_j}{|x_j|^{\varrho-\varepsilon}} = \infty \mbox{ \ws for all \ws } \varepsilon > 0, \end{equation}
then with $\mu, d > 0$ arbitrary
\begin{equation}v(x)= v_{X,M,\mu,d}(x):= e^{d|x|^{\varrho + \mu}}\prod_{j=1}^{\infty}\left|1-\frac{x}{x_j}\right|^{m_j} .\end{equation} \end{defi}

After the definitions of the weights we begin to deal with the description of the functions we need for giving a good choice of points and an omission system.

\remark

In \cite{lelu}, Lemma 5.1. (b) states that
\begin{equation}t^A \leq \frac{tQ^{'}(tx)}{Q^{'}(x)} \leq t^B, \ws \ws x \in (0, \infty), \ws t\in (1, \infty), \end{equation}
and
\begin{equation}A \leq \frac{xQ^{'}(x)}{Q(x)} \leq B, \ws \ws x \in (0,\infty)\end{equation}
Together with the definition it means that on  $(0,\infty)$  $Q^{'} > 0 $;  $Q(cx) \sim Q(x) (c>0)$;  $Q^{'}(x) \sim \frac{Q(x)}{x},$ where $f(x) \sim g(x)$ means that there are positive constants $C$ and $D$ such that $f(x) \leq C g(x)$ and $g(x) \leq D f(x).$ So this is the inspiration of the following definition:

\begin{defi} $f$  grows "polynomially uniformly" if it is three times differentiable and $f^{'}$ is positive and convex on $(0,\infty)$, and  there exists an $x_0>0$ such that on  $(x_0,\infty)$ the followings are valid:
\begin{equation}  f(cx) \sim f(x) \ws \ws (c>0) \end{equation}
\begin{equation}  f^{'}(x) \sim \frac{f(x)}{x} \end{equation} \end{defi}

With this property we can define an admissible function and point system as it follows:

\begin{defi}$w$ is a Freud weight, $Q = \log\frac{1}{w}$, and let us suppose for simplicity that $\frac{Q(x)}{x^3}$ is quasimonotone, that is there exists a monotone function: $m(x)$; for which $m(x)\sim \frac{Q(x)}{x^3}$ on $(x_0, \infty)$ with some $x_0$, furthermore let $\frac{3}{2}< A \leq B$, $\gamma > 0$ is a positive number.$g$ is an admissible function with respect to $Q$ and $\gamma$, if it grows polynomially uniformly on  $(0,\infty)$, and the following relations are valid:
\begin{equation}\frac{g^{[-1]}(x)}{Q^{[-1]}(x)} =O\left(\frac{1}{x^{2\gamma}}\right),  \end{equation}
and there is an $x_0>0$ and an $\varepsilon >0$ such that $\frac{g^{[-1]}(x)}{(Q^{[-1]}(x))^{1-\varepsilon}}$ is decreasing on $(x_0,\infty)$;
\begin{equation}x^{\delta}\max\left\{\frac{(Q^{[-1]}(g(x)))^{\frac{1}{4}}}{(g(x))^{\frac{1}{6}}};\frac{1}{(Q^{[-1]}(g(x)))^{\frac{1}{2}}}\right\}\longrightarrow 0 \ws \ws\mbox{with a} \ws\ws \delta > \frac{5}{4}\end{equation}
when $x \longrightarrow \infty .$

\end{defi}

(Here $g^{[-1]}$ denotes the inverse of the function $g$.) 

\begin{defi}$M$ is an admissible system of positive numbers with respect to $\gamma$, if $ 0 < m_j < 1+\gamma$, and $\liminf_{j \to \infty}m_j >0$.\end{defi}

After these definitions and notations we can formulate the main theorem:

\begin{theorem}Let $w$ be a Freud weight on the real line with the properties were given in Definition 4, and $ 0 < \gamma < \frac{1}{2B}$. Furthermore let $g$ is an admissible function with respect to $Q$ and $\gamma$, and $M$ is an admissible system of positive numbers with respect to $\gamma$. In this situation there is a point system $X \subset \mathbb{R}$
and an "omission system": $\Psi_k = p_{l_k}(w)w$ with
\begin{equation}l_k = g(k) + O(k), \end{equation}
\begin{equation}\mbox{and}\ws\ws\ws \ws\ws\ws\ws\ws\ws\ws\ws\ws\ws\ws\ws d, \mu > 0,\end{equation}
such that the system
\begin{equation}\{\varphi_l\}_{l=1}^{\infty}:= \{p_k(w)w\}_{k=0}^{\infty}\setminus \{\Psi_n\}_{n=1}^{\infty}\end{equation}
is complete and minimal in $L_{v_{X,M,\mu,d}}^p$, where $\inf_{m_j<1}\frac{1}{1-m_j}>p >\sup_{j \atop \gamma - m_j< 0}\frac{1}{\gamma - m_j+1}$, if there are some $m_j$-s for which $\gamma - m_j <0$, and for $\inf_j\frac{1}{1-m_j}> p > 1$, if $\gamma - m_j \geq 0$ for all $j$. \end{theorem}

\remark

With the assumptions of the theorem we will be able to give a numerical method to compute the conjugate system.

\noindent {\bf Examples}

\noindent (1)  Let $Q(x) = |x|^{\beta}; g(x)= x^{\alpha},\alpha, \beta >1.$ It is admissible if  $ \gamma < \frac{\alpha -\beta}{2\alpha \beta}$, and $\alpha > \max\{\frac{15\beta}{2\beta-3};\frac{5\beta}{2}\}$.

eg. for Hermite weight, that is $\beta = 2$, we can choose eg. $\gamma = \frac{1}{5}$ and $\alpha = 31.$ Or if $\beta = 6$,  let $\gamma = \frac{1}{24}$ and $\alpha = 16.$ 

\noindent (2)  Let $Q(x) = |x|^{\beta}, \beta >1$ again, and let $f(x) = x^{\nu}\log x; g(x)= x^{\alpha}\log x.$ 

In this case the relations are the same as in the previous case.

\section{Proof}

As we have seen in the introduction at first we have to solve the following infinite systems of linear equations:

\begin{equation} \left[
\begin{array}{ccccc}
                \Psi_{1}(x_1) & \Psi_2(x_1) & \dots & \Psi_n(x_1)& \dots\\
                \Psi_1(x_2) & \Psi_2(x_2) & \dots & \Psi_n(x_2) & \dots\\
                \vdots &     \vdots &      &  \vdots &      \\
                \Psi_1(x_k) & \Psi_2(x_k) & \dots &  \Psi_n(x_k) & \dots\\
                \vdots &     \vdots &      &  \vdots &     \end{array} \right]
\left[ \begin{array}{c} a_{1m}\\
                a_{2m}\\
                \vdots\\
                a_{km}\\
                \vdots \end{array} \right] = 
\left[ \begin{array}{c} \varphi_m(x_1)\\ 
                \varphi_m(x_2)\\ 
                \vdots\\
                \varphi_m(x_k)\\
                \vdots\end{array}\right] \end{equation}
                
denoted by $Aa_m = c_m$. 

Inconnection with this infinite linear equation system we have to deal with two questions: to get some solution, and to guarantee the convergence of the solution in some sense. Together with the convergence, the exictence of the solution results a biorthonormal system with respect to $\{\varphi_l\}_{l=1}^{\infty}$, and the uniqueness of the solution ensures the completeness of $\{\varphi_l\}_{l=1}^{\infty}$.

\subsection{Solvability}
\subsubsection{Existence}

For the first problem we have to cite the theorem O. Toeplitz \cite{top}, \cite{ban}

\noindent {\it {\bf Theorem A}
The necessary and sufficient condition of the existence a solution of an infinite linear equation system 
$$\sum_{k=1}^{\infty} a_{ki}x_k = y_i,\ws \ws i=1,2,\dots ,$$
is the following: for all $r$ natural, and $h_1, h_2, \dots , h_r$ real numbers for which $\sum_{i=1}^rh_i a_{ki} = 0, \ws \ws k=1,2,\dots ,$ the equality $\sum_{i=1}^rh_i y_i = 0$ fulfils; in particular if the condition $\sum_{i=1}^rh_i a_{ki} = 0, \ws \ws k=1,2,\dots ,$ implies that $h_1=h_2= \dots h_r=0$ the above equation system has a solution for all $\{y_i\}$-s.}

Now we can define our point- and our omission system. For this construction and in the followings we need the following
notion of Mhaskar-Rahmanov-Saff number with respect to $w$, which shows "Where does the sup-norm of a weighted polynomial live" \cite{ms}               
 \begin{defi}$w$ is a Freud weight on the real line. $a_n = a_n(w)$ is the MRS number associated with $w$ if for all $q_n$ polynomials with degree $n$ the followings are valid:
 \begin{equation}\|q_n\|_{w,\infty} = \max_{|x|\leq a_n}|q_n(x)w(x)|,\end{equation}
 and
  \begin{equation}\|q_n\|_{w,\infty} > |q_n(x)w(x)|\ws\ws \mbox{ for  all}\ws \ws |x|>a_n\end{equation}\end{defi}
 
 \remark
 
 $a_u$ is the positive root of the equation
 $$u = \frac{2}{\pi}\int_0^1 a_utQ^{'}(a_ut)(1-t^2)^{-\frac{1}{2}}dt, \ws\ws u>0.$$

 \begin{lemma} Let $Q$ and $\gamma$ be as in Theorem 1, and let $g$ is an admissible function with respect to $Q$ and $\gamma$. Now there is a point system $X \subset \mathbb{R}$, and an omission system $\Psi_k = p_{l_k}(w)w$ with
$l_k = g(k) + O(k)$, such that 
\begin{equation}|\Psi_{i-1}(x_{i})| > c \|\Psi_{i}\|_{\infty} \ws\ws\ws i=1,2,\dots\end{equation}
with an absolute constant $c$, and the determinants
\begin{equation} D_n = \left|
                \begin{array}{cccc}
                \Psi_{1}(x_1) &  \Psi_2(x_1) &\dots & \Psi_n(x_1)\\
                \Psi_1(x_2) & \Psi_2(x_2) & \dots & \Psi_n(x_2) \\
                \vdots &     \vdots &      &  \vdots      \\
                \Psi_1(x_n) & \Psi_2(x_n) & \dots &  \Psi_n(x_n) \\
                \end{array} \right| \neq 0 \ws \ws n \in\mathbb{N} \end{equation} \end{lemma}

\proof

Let $x_1 \in \bf{R}_+$ be an arbitrary point, say $x_1=1,$ and let $n_0 \in \mathbb{N}$ be a fixed number (will be given later). Let $g^{+}$ be a function with the properties of $g$, and let us denote by 
$$g(k)= g^{*}(n_0+k).$$
We can choose $\Psi_1 = p_{k_1}w$ such that $k_1 = g(1)+ O(1) $ and $\Psi_1(x_1) \neq 0.$
Now let us suppose that $x_1, \dots ,x_n$ and $\Psi_1, \dots , \Psi_n$ had already chosen, such that and $l_k = g(k) + O(k)$, $|\Psi_{k-1}(x_{k})| = \|\Psi_{k}\|_{\infty}$ for $k = 2, \dots ,n$ and $D_k \neq 0$ for $k=1,\dots , n$.

At first we will give $x_{n+1}$ such that $|\Psi_{n}(x_{n+1})| = \|\Psi_{n}\|_{\infty}$. 
So with this choice we get a not too small element in every rows. 
It follows from \cite{lelu} Lemma 5.1 that 
\begin{equation} a_n \sim Q^{[-1]}(n), \end{equation}
and so by the assumptions on $g$ we get that
\begin{equation} |x_k| \sim Q^{[-1]}(g(k)), \end{equation}
 
In the followings we will show that for every $m > l_n$ among the indices $m, m+1, \dots , m+2n+1$ we can find a "good" one, that is there is a $k \in \{m, m+1, \dots , m+2n+1\}$ like that if we choose $\Psi_{n+1} = p_kw$, then $D_{n+1} \neq 0.$ By (8) it means that we can chose $\Psi_{n+1}= p_{l_{n+1}}w$ such that $l_{n+1}= g(n+1) + O(n+1).$

So let us suppose indirectly, that there is an $m > l_n$ for which
\begin{equation} D_{n+1} = \left|
                \begin{array}{cccccc}
                \Psi_{1}(x_1) &  \Psi_2(x_1) &\dots & \Psi_n(x_1) & p_kw(x_1)\\
                \Psi_1(x_2) & \Psi_2(x_2) & \dots & \Psi_n(x_2) & p_kw(x_2)\\
                \vdots &     \vdots &       &    \vdots            & \vdots      \\
                \Psi_1(x_n) & \Psi_2(x_n) & \dots &  \Psi_n(x_n) & p_kw(x_n) \\
                \Psi_1(x_{n+1}) & \Psi_2(x_{n+1}) & \dots &  \Psi_n(x_{n+1}) & p_kw(x_{n+1})
                \end{array} \right| =0  \end{equation}
for all $k \in \{m, m+1, \dots , m+2n+1\}$ 
Let us expand this determinant by the elements of the last column:
\begin{equation}D_{n+1} = (-1)^{n+1}\sum_{j=1}^{n+1}(p_kw)(x_j)(-1)^jB_j = 0,\end{equation}
where $B_j$ is that subdeterminant which comes when the last column and the $j^{th}$ row are omitted. ($B_{n+1} = D_n$.)
Denoting by $A_j := (-1)^j\hat{B}_j $, where  $\hat{B}_j$ are the determinants $B_j$ divided by the product of $w(x_i)$-s we get that 
\begin{equation}\sum_{j=1}^{n+1}p_k(x_j)A_j =0, \ws k \in \{m, m+1, \dots , m+2n+1\}.\end{equation}
Let us recall the recurrence formula of the orthonormal polynomials with respect to the even weight $w$:
\begin{equation} xp_{n+1} =  \varrho_{n+2}p_{n+2}+ \varrho_{n+1}p_{n},\end{equation}
where $\varrho_{n} \sim a_n$ are constants. By this formula we get from (24) that for any $0 \leq l \leq 2n-1$
\begin{equation}\sum_{j=1}^{n+1}x_jp_{m+l+1}(x_j)A_j = \varrho_{m+l+2}\sum_{j=1}^{n+1}p_{m+l+2}(x_j)A_j + \varrho_{m+l+1}\sum_{j=1}^{n+1}p_{m+l}(x_j)A_j = 0 \end{equation}
 by the same argument we have that
 \begin{equation} c_p\sum_{j=1}^{n+1}x_j^p p_{m+l+p}(x_j)A_j = 0 ,\ws \ws 0\leq p \leq n, \ws\ws 0\leq l \leq 2n+1-2p \end{equation}
 that is
 \begin{equation} 0 =\sum_{p=0}^{n} c_p\sum_{j=1}^{n+1}x_j^p p_{k}(x_j)A_j = \sum_{j=1}^{n+1}p_{k}(x_j)A_j\sum_{p=0}^{n} c_px_j^p  = \sum_{j=1}^{n+1}q_n(x_j)p_{k}(x_j)A_j, \end{equation}
 where  $ k=m+n, m+n+1$, and $q_n$ is a polynomial with degree $n$. So let us choose $q_n= q_{n,k}$ like
 \begin{equation} {\mbox sign} \ws q_{n,k}(x_j) = {\mbox sign} \ws p_{k}(x_j)A_j \end{equation}
 (If in a point $x_j$ the expression $p_{k}(x_j)A_j$ is zero, then we have no assumption on the sign of $q_{n,k}$ at $x_j$.) With this choice we get that all the terms of the above sum are zero, but we know that $A_{n+1}=(-1)^{n+1}\hat{D}_n \neq 0$ and we can suppose that $q_{n,k}(x_{n+1})\neq 0$, that is $p_{k}(x_{n+1})$ must be zero for $ k=m+n, m+n+1$, which is impossible, because two consecutive orthogonal polynomials can't have zero at the same point. So the first lemma is proved.
 
 \vspace{3mm}
 
 \note
 
 Denoting by $\Psi_{0}(x_1):=\Psi_{1}(x_1)$ we can define the modified linear equation systems, which are equvivalent with the original ones: 
 \begin{equation} \hat{A}a_m = \hat{c}_m:\end{equation}
 
 \begin{equation} \left[
 \begin{array}{ccccc}
                 \frac{\Psi_{1}(x_1)}{\Psi_{0}(x_1)} & \frac{\Psi_2(x_1)}{\Psi_{0}(x_1)} & \dots & \frac{\Psi_n(x_1)}{\Psi_{0}(x_1)}& \dots\\
                 \frac{\Psi_1(x_2)}{\Psi_{1}(x_2)} & \frac{\Psi_2(x_2)}{\Psi_{1}(x_2)} & \dots & \frac{\Psi_n(x_2)}{\Psi_{1}(x_2)} & \dots\\
                 \vdots &     \vdots &      &  \vdots &      \\
                 \frac{\Psi_1(x_k)}{\Psi_{k-1}(x_k)} & \frac{\Psi_2(x_k)}{\Psi_{k-1}(x_k)} & \dots &  \frac{\Psi_n(x_k)}{\Psi_{k-1}(x_k)} & \dots\\
                 \vdots &     \vdots &      &  \vdots &     \end{array} \right]
 \left[ \begin{array}{c} a_{1m}\\
                 a_{2m}\\
                 \vdots\\
                 a_{km}\\
                 \vdots \end{array} \right] = 
 \left[ \begin{array}{c} \frac{\varphi_m(x_1)}{\Psi_{0}(x_1)}\\ 
                 \frac{\varphi_m(x_2)}{\Psi_{1}(x_2)}\\ 
                 \vdots\\
                 \frac{\varphi_m(x_k)}{\Psi_{k-1}(x_k)}\\
                 \vdots\end{array}\right], \end{equation}

 And let us denote the elements of $\hat{A}\hat{A}^{T}$ by
 \begin{equation} \alpha_{ij} =<{}_i\hat{A}, {}_j\hat{A}> = \sum_{k=1}^{\infty}\frac{\Psi_k(x_i)}{\Psi_{i-1}(x_i)}\frac{\Psi_k(x_j)}{\Psi_{j-1}(x_j)}, \end{equation}
 where $<\cdot,\cdot>$ denotes the usual inner product, and ${}_i\hat{A}$ is the $i^{th}$ row of $\hat{A}$,
 and by $B^{(n)}$ the principal minor of $\hat{A}\hat{A}^{T}$:
 \begin{equation}  B^{(n)} =\left[
 \begin{array}{ccc}\alpha_{11}& \dots & \alpha_{1n}\\
                      \vdots &     &\vdots \\
                       \alpha_{11}& \dots & \alpha_{1n} \end{array}\right], \end{equation} 
 
 and by
 \begin{equation}  \hat{c}_m^{(n)} =\left[
 \begin{array}{c}\frac{\varphi_m(x_1)}{\Psi_{0}(x_1)}\\ 
                 \frac{\varphi_m(x_2)}{\Psi_{1}(x_2)}\\ 
                 \vdots\\
                 \frac{\varphi_m(x_n)}{\Psi_{n-1}(x_n)}
                \end{array}\right]. \end{equation} 
                
 With these notatios we are in the position to formulate the theorem of F. Riesz \cite{riesz}, which will be our basic tool for proving some convergence property of the solution.

 \noindent {\it {\bf Theorem B} 
 
 With the notation
 \begin{equation} M^*(\frac{\varphi_m(x_1)}{\Psi_{0}(x_1)},\frac{\varphi_m(x_2)}{\Psi_{1}(x_2)}, \dots ) = \lim_{n \to \infty}\left( - \frac{\left|\begin{array}{cc}B^{(n)} & \hat{c}_m^{(n)}\\
                                                   \left(\hat{c}_m^{(n)}\right)^T & 0\end{array}\right|}
                                                  {| B^{(n)}|}\right)^{\frac{1}{2}} \end{equation}
The equation $\hat{A}a_m = \hat{c}_m$ has a solution for which
\begin{equation} \|a_m\|_2 = \left(\sum_{k=1}^{\infty}a_{km}^2 \right)^{\frac{1}{2}} \leq M \end{equation}
iff
 \begin{equation} M^*(\frac{\varphi_m(x_1)}{\Psi_{0}(x_1)},\frac{\varphi_m(x_2)}{\Psi_{1}(x_2)}, \dots ) \leq M \end{equation}  }

 For an estimation on $M^*$ we need some lemmas. At first we have to introduce a 
 
 \note
 
 Let us define a function
 \begin{equation}f(c_0,\delta) = \frac{(6c_0^2+2c_0)^2}{4^{\delta}}\left[1+\frac{2\sqrt{2}}{\delta-\frac{5}{4}}\exp\left((6c_0^2+2c_0)^2\frac{15}{6}\right)\right],\end{equation}
 where $\delta > \frac{5}{4}$ arbitrary, and let us give $c_0 = c_0(\delta)$ such that $f(c_0,\delta) \in (0,1).$ 
 
 Furthermore let us denote by $w_i= w_i^{(n)}$ the $i^{th}$ row of the symmetric matrix $B^{(n)}$. 
 
 \begin{lemma}  
 With the previous notations,
 if there is a $\delta >\frac{5}{4}$ and a $c<c_0(\delta)$, such that
 \begin{equation} |\alpha_{ij}| < c \max\{i,j\}^{-\delta}, \ws\ws\ws \mbox{if} \ws \ws \ws i\neq j\end{equation}
  then there exists a $0 < q <1$ for which 
  \begin{equation} \left|B^{(n)}\right| > q \prod_{j=1}^n \|w_j\|_2. \end{equation} \end{lemma}
  
  \proof
  
  Let us suppose that $i<j$. At first we will prove that assuming that $|\alpha_{ij}| < c j^{-\delta},$  the  cosine of the angle of the $i^{th}$ and $j^{th}$ rows is of order $ j^{-\delta}$, that is
  \begin{equation} \left|\cos \beta_{ij}\right| :=\frac{|<w_i,w_j>|}{\|w_i\|_2\|w_j\|_2}  < c_1 j^{-\delta},\end{equation}
  where $c_1 = 6c^2+2c$. Let us observe at first that  $\|w_j\|_2 \geq \alpha_{jj}= \|{}_j\hat{A}\|_2^2 \geq 1$. Thus we get that
  
  $$\frac{|<w_i,w_j>|}{\|w_i\|_2\|w_j\|_2} \leq \sum_{k=1}^{i-1}\left|<{}_i\hat{A},{}_k\hat{A}><{}_j\hat{A},{}_k\hat{A}>\right| + \sum_{k=i+1}^{j-1}|\cdot | + \sum_{k=j+1}^n|\cdot|$$ 
  $$+ \left|<{}_i\hat{A},{}_j\hat{A}>\right|\frac{\|{}_i\hat{A}\|_2^2 + \|{}_i\hat{A}\|_2^2}{\|w_i\|_2\|w_j\|_2} $$
  \begin{equation}\leq c^2 i^{1-\delta}j^{-\delta} + c^2 j^{-\delta}\frac{i^{1-\delta}}{\delta-1} + c^2 \frac{j^{1-2\delta} }{2\delta -1} + cj^{-\delta}\frac{\|{}_i\hat{A}\|_2^2 + \|{}_i\hat{A}\|_2^2}{\|{}_i\hat{A}\|_2^2\|{}_j\hat{A}\|_2^2}  \leq c_1 j^{-\delta}\end{equation}
 In the last step we used that $\delta > 1.$
 
 By this inequality we can prove the original one. With the notation
 $$B^{(n)}_0 = \left[\begin{array}{cccc}\frac{w_1^T}{\|w_1\|_2},&\frac{w_2^T}{\|w_2\|_2},&\dots ,& \frac{w_n^T}{\|w_n\|_2}\end{array}\right],$$
 we have to show that $\left|B^{(n)}_0\right| \geq q$. Let us estimate
 $$\left|\left(B^{(n)}_0\right)^TB^{(n)}_0\right|=\left|\begin{array}{cccccc}1 &  & \dots & & &\cos\beta_{1n}\\
                                                                            \vdots & \ddots &  & \cos \beta_{ij} & &\vdots \\
                                                                            \dots &  &  &  & &\dots \\
                                                                            \vdots &  & & \ddots  & & \vdots \\                                                                                                                                                                                                                                  \cos\beta_{1n}&  &  & & &1 \end{array}\right|
   $$ $$=\sum_{k=1}^{n-1}(-1)^{k+n}\cos \beta_{kn}\det B_k^{(n)} + \det ((B_0^{(n-1)})^TB_0^{(n-1)}),$$
   where
   $$B_k^{(n)}=\left[\begin{array}{cccccc}1 &\cos\beta_{12}  & \dots & & &\cos\beta_{1n-1}\\
                                                                            \vdots & \ddots &  & \cos \beta_{ij} & &\vdots \\
                                                                            \cos\beta_{k-11} &  &  &  & &\cos\beta_{k-1n-1} \\
                                                                            \cos\beta_{k+11} &  &  &  & &\cos\beta_{k+1n-1} \\
                                                                            \vdots &  & & \ddots  & & \vdots \\                                                                                                                                                                                                                                  \cos\beta_{n1}&  &  &  \dots  & &\cos\beta_{nn-1} \end{array}\right]$$
             
 by Hadamard's inequality we get that
 $$\left|\det B_k^{(n)}\right| \leq \prod_{i=1\atop i\neq k}^{n-1}\sqrt{1+\sum_{l=1\atop l\neq i}^{n-1}\cos^2 \beta_{il}} \sqrt{\sum_{l=1}^{n-1}\cos^2 \beta_{nl}}  
 $$ $$\leq \prod_{i=1\atop i\neq k}^{n-1}\sqrt{1+(i-1)c_1^2i^{-2\delta} + \sum_{k=1}^{n-1-i}\frac{1}{(i+k)^{2\delta}}} \sqrt{c_1^2\frac{n-1}{(n)^{2\delta}}}$$ $$ \leq c_1 \frac{\sqrt{n-1}}{(n)^{\delta}}\prod_{i=1\atop i\neq k}^{n-1}\sqrt{1+c_1^2i^{1-2\delta}\left(1+\frac{1}{2\delta - 1}\right)} \leq c_1 n^{\frac{1}{2}-\delta}\prod_{i=1\atop i\neq k}^{n}\sqrt{1+\frac{5}{3}c_1^2i^{1-2\delta}}$$ $$\leq c_1 n^{\frac{1}{2}-\delta}\exp\left(\frac{5}{6}c_1^2\left(1+\frac{1}{2(\delta - 1)}\right)\right)$$
 \begin{equation}\leq c_1 \exp\left(c_1^2\frac{15}{6}\right)n^{\frac{1}{2}-\delta}= c_2n^{\frac{1}{2}-\delta}\end{equation}
  By this calculation we get that
  $$\det ((B_0^{(n)})^TB_0^{(n)}) \geq \det ((B_0^{(n-1)})^TB_0^{(n-1)}) - (n-1)c_1c_2n^{-\delta}n^{\frac{1}{2}-\delta} $$ $$\geq \det ((B_0^{(n-2)})^TB_0^{(n-2)}) - c_1c_2(n-1)^{\frac{3}{2}-2\delta}-c_1c_2n^{\frac{3}{2}-2\delta} $$ $$\geq \dots \geq \det ((B_0^{(2)})^TB_0^{(2)}) -c_1c_2\sum_{k=3}^n k^{\frac{3}{2}-2\delta}$$ \begin{equation}\geq \det ((B_0^{(2)})^TB_0^{(2)}) -c_1c_2\frac{1}{2^{2\delta-\frac{5}{2}}\left(2\delta-\frac{5}{2}\right)} \geq 1 - \frac{c_1^2}{4^{\delta}}-c_1c_2\frac{1}{2^{2\delta-\frac{5}{2}}\left(2\delta-\frac{5}{2}\right)}\end{equation}
  By the last inequality, the assumptions on $c_0$ implies that there is a $q_1 \in (0,1)$ such that $(\det B_0^{(n)})^2 >q_1,$ which proves the lemma.
   
 \corollary
 
 With the assumptions of Lemma 2 and the notations above, the following inequality is valid:
 \begin{equation}- \frac{\left|\begin{array}{cc}B^{(n)} & \hat{c}_m^{(n)}\\
                                                   \left(\hat{c}_m^{(n)}\right)^T & 0\end{array}\right|}
                                                  {| B^{(n)}|}
    \leq \frac {1}{q} \|\hat{c}_m\|_2 e^{\frac{1}{2}\|\hat{c}_m\|_2^2} \end{equation}  
 \proof
 
 Applying Hadamard's inequality again, and recalling that $\|w_i\|_2 \geq 1$, we get that
 $$- \frac{\left|\begin{array}{cc}B^{(n)} & \hat{c}_m^{(n)}\\
                                                   \left(\hat{c}_m^{(n)}\right)^T & 0\end{array}\right|}
                                                  {| B^{(n)}|}
    \leq \frac{1}{q}\frac{\prod_{j=1}^n \sqrt{\sum_{i=1}^n \alpha_{ij}^2 + \hat{c}_{m,j}^2}\sqrt{\sum_{i=1}^n\hat{c}_{m,i}^2}}{\prod_{j=1}^n\sqrt{\sum_{i=1}^n\alpha_{ij}^2}} $$
    $$\leq \frac {1}{q}\|\hat{c}_m\|_2\left(\prod_{j=1}^n \left(1+ \frac{\hat{c}_{m,j}^2}{\|w_j\|_2^2}\right)\right)^{\frac{1}{2}} \leq \frac {1}{q}\|\hat{c}_m\|_2 e^{\frac{1}{2}\sum_{j=1}^n\hat{c}_{m,j}^2},$$
    and the corollary is proved.
 
 \begin{lemma}
 With the previous notations, there is a $\delta > \frac{5}{4}$ and a $c_0 = c_0(\delta)$ such that $f(c_0,\delta) \in (0,1)$ (see (38) for $f(c_0,\delta)$), with which
 \begin{equation}|\alpha_{ij}| \leq c\max\{i, j\}^{-\delta}, \ws\ws\ws \mbox{if} \ws \ws \ws i\neq j\end{equation}
for a $c<c_0$ .\end{lemma}

 \remark

  \noindent (A)We can easily deduce eg. from 2.19 of \cite{lelu2} or 2.6 of \cite{mhas}, that 
  \begin{equation}|p_n(w,x)w(x)| \leq c \sqrt{w(x)} \ws \ws {\mbox if} \ws \ws |x| > (1+c)a_n .\end{equation}
  
  \noindent (B)\cite{lelu} Cor. 1.4 : If $w$ is a Freud weight, then
  \begin{equation}\sup_{x\in \mathbb{R}}|p_n(w,x)|w(x)\left|1 - \frac{|x|}{a_n}\right|^{\frac{1}{4}} \sim a_n^{-\frac{1}{2}} \end{equation}
  and
  \begin{equation}\sup_{x\in \mathbb{R}}|p_n(w,x)|w(x) \sim n^{\frac{1}{6}}a_n^{-\frac{1}{2}}\end{equation}

 \note
  
  Let us denote by $I_{MRS}(p_k)$ the support of the equlibrium measure with respect to $w^k$, that is
 \begin{equation}I_{MRS}(p_k)= I_{MRS}(p_kw) = [-a_k,a_k],\end{equation}
 
 \proof
 
 $$|\alpha_{ij}| \leq \frac{1}{|\Psi_{i-1}(x_i)||\Psi_{j-1}(x_j)|}$$
 $$\times \left(\sum_{k=1}^{\frac{1}{c}j}|\Psi_k(x_i)\Psi_k(x_j)| + \sum_{k=\frac{1}{c}j+1}^{j-2}| \cdot | + \sum_{k=j-1}^{cj}| \cdot | + \sum_{k=cj+1}^{\infty}| \cdot |\right) $$
 \begin{equation}= \frac{1}{|\Psi_{i-1}(x_i)||\Psi_{j-1}(x_j)|}(S_1 + S_2 + S_3 + S_4)\end{equation}
 At first, recalling the special assumption on the denominator, considering (48) we have that
 $$ |\Psi_{i-1}(x_i)| \sim  \|\Psi_{i-1}\|_{\infty} =  \|p_{l_{i-1}}(w)\|_{\infty} \sim \left(l_{i-1}\right)^{\frac{1}{6}}a_{l_{i-1}}^{-\frac{1}{2}}$$
 \begin{equation} \sim g(i)^{\frac{1}{6}}\left(Q^{[-1]}(g(i))\right)^{-\frac{1}{2}}\end{equation}
 
 By the first remark, the members in $S_1$ are exponentially small, because eider $x_i$ and $x_j$, or only $x_j$, are out of $cI_{MRS}(\Psi_k)$ for such $k$-s. According to the previous calculation we get that
 $$\frac{S_1}{|\Psi_{i-1}(x_i)||\Psi_{j-1}(x_j)|}$$ $$ \leq c e^{-\frac{Q}{2}(x_j)}j(g(i))^{-\frac{1}{6}}(Q^{[-1]}(g(i)))^{\frac{1}{2}}(g(j))^{-\frac{1}{6}}(Q^{[-1]}(g(j)))^{\frac{1}{2}}=(*)$$
 Now we have to distinguish two cases according to the infinite norm of the weighed orthonormal polynomials tend to infinity with the degree of the polynomial, or it is bounded (see \cite{lelu} 5.1, and the assumption in Definition 4). That is in the second case
 \begin{equation}(*)\leq c e^{-cg(j))}j(g(j))^{-\frac{1}{6}}(Q^{[-1]}(g(j)))^{\frac{1}{2}},\end{equation}
 and in the first case
 \begin{equation}(*)\leq c e^{-cg(j))}j(g(j))^{-\frac{1}{3}}(Q^{[-1]}g(j)).\end{equation}
 (Here all the $c$-s are different absolute constants.) So in both cases
  \begin{equation}\frac{S_1}{|\Psi_{i-1}(x_i)||\Psi_{j-1}(x_j)|}\leq cj^{-\delta} \end{equation}
 for a $\delta >\frac{5}{4}$ with a $c <c_0$, if and only if 
 $$c je^{-cg^{*}(n_0+j))}$$
 $$ \times\max\{(g^{*}(n_0+j))^{-\frac{1}{6}}(Q^{[-1]}(g^{*}(n_0+j)))^{\frac{1}{2}},(g^{*}(n_0+j))^{-\frac{1}{3}}(Q^{[-1]}g^{*}(n_0+j))\}$$
 \begin{equation}\leq \frac{1}{4} c_0 j^{-\delta}\end{equation}
 That is we can choose an $n_0$ enough large with which (46) will be valid.

 In $S_4$ we collected that terms which maximum points are far away from $x_i$ and $x_j$. Applying \cite{lelu} (1.20):
 $$S_4 \leq c \sum_{k=cj+1}^{\infty}\frac{1}{\sqrt{a_{l_k}}}\frac{1}{\left(a_{l_k}-|x_i|\right)^{\frac{1}{4}}\left(a_{l_k}-|x_j|\right)^{\frac{1}{4}}} \leq c \sum_{k=cj+1}^{\infty}\frac{1}{a_{l_k}}$$
 $$\leq c \int_{cj+1}^{\infty}\frac{1}{Q^{[-1]}(g(x))}dx \leq c \int_{cQ^{[-1]}(g(j))}^{\infty}\frac{\left(g^{[-1]}(Q(y))\right)^{'}}{y}dy = (*)$$
 Using the polynomially growing property of $g$ and $Q$, and then the monotonicity of the lefthand side of (10), we get that
 $$(*) \leq c \int_{cQ^{[-1]}(g(j))}^{\infty}\frac{g^{[-1]}(Q(y))}{y^{1-\varepsilon +1+\varepsilon}}dy $$
 $$\leq c \frac{n_0 +j}{(Q^{[-1]}(g^{*}(n_0+j))))^{1-\varepsilon}}\int_{cQ^{[-1]}(g^{*}(n_0+j))}^{\infty}\frac{1}{y^{1+\varepsilon}}dy $$ \begin{equation}\leq c \frac{n_0 +j}{Q^{[-1]}(g^{*}(n_0+j))} \end{equation}
  It means that
  $$ \frac{S_4}{|\Psi_{i-1}(x_i)||\Psi_{j-1}(x_j)|} $$
  $$\leq c \frac{n_0 +j}{Q^{[-1]}(g^{*}(n_0+j))}$$
  \begin{equation}\times  \max\{(g^{*}(n_0+j))^{-\frac{1}{6}}(Q^{[-1]}(g^{*}(n_0+j)))^{\frac{1}{2}},(g^{*}(n_0+j))^{-\frac{1}{3}}(Q^{[-1]}g^{*}(n_0+j))\} \end{equation}
So it is clear that by the assuptions on $g$ and $Q$, and by (11), if $n_0$ is enough large, then
 \begin{equation} \frac{S_4}{|\Psi_{\mu_i}(x_i)||\Psi_{\mu_j}(x_j)|} \leq \frac{1}{4} c_0 j^{-\delta}. \end{equation}
 
 In $S_2$ $x_j \notin I_{MRS}(\Psi_k)$, in $S_3$ $x_i$ and $x_j$ are both in $I_{MRS}(\Psi_k)$, but we can handle the two sums similarly: there are $O(j)$ terms in that sums, and at most one term, the $(i-1)^{th}$ or the $(j-1)^{th}$, has a factor $1$. Furthermore the distance between two consecutive maxima is more than some constant $c$ (Def. 4):
 \begin{equation} a_{l_{k+1}}-a_{l_k}\sim \frac{Q^{[-1]}(g(k))}{k} \end{equation}
 This implies that if $k \neq j-1$, say, we can estimate by (48)
 \begin{equation}|\Psi_{k}(x_j)|\leq c a_{l_k}^{-\frac{1}{4}}\left(a_{l_{k+1}}-a_{l_k}\right)^{-\frac{1}{4}} \leq c j^{\frac{1}{4}}(Q^{[-1]}(g(j)))^{-\frac{1}{2}}  \end{equation}
 Let us assume at first that $\frac{1}{c}j+1 \leq i-1 \leq j-2$, or $i \sim j$ and so
 $$\frac{S_2}{|\Psi_{i-1}(x_i)||\Psi_{j-1}(x_j)|} \leq c \frac{|\Psi_{i-1}(x_j)|}{|\Psi_{j-1}(x_j)|} 
  + \sum_{k=\frac{1}{c}j+1 \atop k\neq i-1}^{j-2}\frac{|\Psi_k(x_i)\Psi_k(x_j)|}{|\Psi_{i-1}(x_i)||\Psi_{j-1}(x_j)|} $$ \begin{equation}\leq c \frac{j^{\frac{1}{4}}}{(g(j))^{\frac{1}{6}}}  + cj\frac{j^{\frac{1}{2}}}{(g(j))^{\frac{1}{3}}}=(*),\end{equation}
  and again by the assuptions on $g$ and $Q$ (11), if $n_0$ is enough large, then
  \begin{equation} (*) \leq \frac{1}{4} c_0 j^{-\delta}. \end{equation}
 
 If $i<<j$, then the first term is missing, and $a_{l_k}-a_{l_{i-1}}>ca_{l_k}$.
 $$\frac{S_2}{|\Psi_{i-1}(x_i)||\Psi_{j-1}(x_j)|} \leq c \frac{(Q^{[-1]}(g(i)))^{\frac{1}{2}}}{(g(i))^{\frac{1}{6}}} \frac{(Q^{[-1]}(g(j)))^{\frac{1}{2}}}{(g(j))^{\frac{1}{6}}}j \frac{j^{\frac{1}{4}}}{Q^{[-1]}(g(j))}$$
 If the first member is bounded, then 
 $$\frac{S_2}{|\Psi_{i-1}(x_i)||\Psi_{j-1}(x_j)|} \leq c \frac{j^{\frac{1}{4}}}{(g(j))^{\frac{1}{6}}}\frac{j}{(Q^{[-1]}(g(j)))^{\frac{1}{2}}}$$ 
 If it can be estimated by the second, then
 $$\frac{S_2}{|\Psi_{i-1}(x_i)||\Psi_{j-1}(x_j)|} \leq c \frac{j^{\frac{5}{4}}}{(g(j))^{\frac{1}{3}}}$$
 so according to (11)
 \begin{equation}\frac{S_2}{|\Psi_{i-1}(x_i)||\Psi_{j-1}(x_j)|} \leq  \frac{1}{4}c_0 j^{-\delta} ,\end{equation}
 if $n_0$ is enough large.
 
 We can estimate $S_3$ on the same way. Here the exceptional term is
 $$\frac{|\Psi_{j-1}(x_i)|}{|\Psi_{i-1}(x_i)|}.$$
 If $\|\Psi_{i-1}\|_{\infty} \geq c \|\Psi_{j-1}\|_{\infty}$ (eg. if $i\sim j$ or if the infinite norm of the weighed orthonormal polynomials tend to zero), then the above term can be estimated by $c \frac{j^{\frac{1}{4}}}{(g(j))^{\frac{1}{6}}}$ as in $S_2$. 
 
 If $i<<j$, and the reciprocal of the infinite norm of the weighed orthonormal polynomials is bounded, then $a_{l_{j-1}}-a_{l_{i-1}} > c a_{l_{j-1}}$, and by (52)
 $$\frac{|\Psi_{j-1}(x_i)|}{|\Psi_{i-1}(x_i)|} \leq c \frac{(Q^{[-1]}(g(i)))^{\frac{1}{2}}}{(g(i))^{\frac{1}{6}}}\frac{1}{(a_{l_{j-1}}(a_{l_{j-1}}-a_{l_{i-1}} ))^{\frac{1}{4}}} $$ \begin{equation}\leq c \frac{1}{(Q^{[-1]}(g(j)))^{\frac{1}{2}}}\leq \frac{1}{4} c_0 j^{-\delta} ,\end{equation}
by (11), if $n_0$ is enough large.
 
 The sum, without the extremal term can be estimated as in $S_2$, and so the lemma is proved.
 
 In the following lemma we state that the operator $\hat{A}$ acts, and is bounded on $l_2$, and $Ran\hat{A} = l_2$.
 
 \begin{lemma}
 With the previous notations for all $\hat{c}_m \in l_2$ there is an $a_m \in l_2$, such that $\hat{A}a_m=\hat{c}_m$, and
 \begin{equation}\|a_m\|_2 \leq c\sqrt{\|\hat{c}_m\|_2}e^{c\|\hat{c}_m\|_2^2}, \end{equation}
 and if $\hat{A}a_m=\hat{c}_m$ with some $a_m \in l_2$, then
 \begin{equation}\|\hat{c}_m\|_2 \leq c\|a_m\|_2 , \end{equation}
 where $c$-s  are different absolute constants. \end{lemma}
 
 \proof
 
 Theorem B, Lemma 2, the corollary after Lemma 2 and Lemma 3 prove (66). For proving (67) let us consider
 \begin{equation}\|\hat{c}_m\|_2 = \left(\sum_{i=1}^{\infty}\hat{c}_{im }^2\right)^{\frac{1}{2}} =
 \left(\sum_{i=1}^{\infty} \left(\sum_{k=1}^{\infty}\frac{\Psi_k(x_i)}{\Psi_{i-1}(x_i)}a_{km}\right)^2\right)^{\frac{1}{2}}.\end{equation}
 Let us decompose the vector $\hat{c}_m$ to two parts:
 $$\hat{c}_m = \hat{c}_m^{(1)}+\hat{c}_m^{(2)},$$
 where
 $$\hat{c}_{im}^{(1)}=\sum_{1\leq k<\infty \atop k \neq i-1}\frac{\Psi_k(x_i)}{\Psi_{i-1}(x_i)}a_{km},$$
 and 
 $$\hat{c}_{im}^{(2)} = a_{i-1m}.$$
It is clear that
 \begin{equation}\|\hat{c}_m^{(2)}\|_2 \leq \|a_m\|_2 .\end{equation}
 According to (68)
 $$\|\hat{c}_m^{(1)}\|_2 \leq \|a_m\|_2 \left(\sum_{i=1}^{\infty} \sum_{1\leq k<\infty \atop k \neq i-1}\left(\frac{\Psi_k(x_i)}{\Psi_{i-1}(x_i)}\right)^2\right)^{\frac{1}{2}}$$
 \begin{equation} \leq c \|a_m\|_2 \left(\sum_{i=1}^{\infty} \frac{(Q^{[-1]}(g(i)))}{(g(i))^{\frac{1}{3}}}\sum_{1\leq k<\infty \atop k \neq i-1}\Psi_k^2(x_i)\right)^{\frac{1}{2}}\end{equation}
 We can decompose the inner sum to three parts:
 $$\sum_{1\leq k<\infty \atop k \neq i-1}\Psi_k^2(x_i)=\sum_{1\leq k <\frac{1}{c}i}\Psi_k^2(x_i) + \sum_{\frac{1}{c}i\leq k< ci \atop k \neq i-1}(\cdot)+ \sum_{ci \leq k<\infty}(\cdot) = S_1+S_2+S_3 $$
 As we have shown in Lemma 3, $S_1$ is exponentially small. Also as in Lemma 3 (in the estimation of $S_4$)
 \begin{equation}S_3 \leq c \sum_{ci \leq k<\infty}\frac{1}{a_{l_k}} \leq c \frac{n_0 +i}{Q^{[-1]}(g^{*}(n_0+i))},\end{equation}
 and as in (61)
 \begin{equation} S_2 \leq c i(Q^{[-1]}(g(i)))^{-\frac{1}{2}}.\end{equation}
 So according to the previous calculation and (11) we obtain that
 \begin{equation} \|\hat{c}_m^{(1)}\|_2 \leq c \|a_m\|_2 \left(\sum_{i=1}^{\infty} \frac{(Q^{[-1]}(g(i)))}{(g(i))^{\frac{1}{3}}} i(Q^{[-1]}(g(i)))^{-\frac{1}{2}}\right)^{\frac{1}{2}}\leq c \|a_m\|_2,\end{equation}
 which proves the lemma. 
 
 \remark
 
 It is wellknown (see eg \cite{mate}), that if $T: H_1 \longrightarrow H_2$ is a continuous linear operator between two Hilbert spaces, then $TT^{*}$ has an inverse, iff $Ran T= H_2$, and in this situation $T^{*}(TT^{*})^{-1}y$ gives the solution with the minimal norm of the linear equation $Tx=y$. Hence we get the following
 
 \corollary
 
 $\hat{A}x=\hat{c}_m$ has a solution $a_m$ in $l_2$ with the minimal norm (and it is unique with this property), and
 \begin{equation} a_m = \hat{A}^T(\hat{A}\hat{A}^T)^{-1}\hat{c}_m\end{equation}
 
 \subsubsection{Unicity}
 
 On the same chain of ideas, by changing the role of $\hat{A}$ and $\hat{A}^T$, we will prove that $\hat{A}^T\hat{A}$ has an inverse on $l_2$, that is {\it Ker}$\hat{A} = \{0\}$ (see eg \cite{mate}). For this we need the following notations and Lemma:
 
 \note
 
 Let us denote by $\lambda_{kl}$ the elements of the matrix $\hat{A}^T\hat{A}$:
 \begin{equation}{}_k\hat{A}^T\hat{A}_l = \lambda_{kl} = \sum_{m=1}^{\infty}\frac{\Psi_k(x_m)\Psi_l(x_m)}{\Psi_{m-1}^2(x_m)} \end{equation}
 
 \remark
 
 As in the previous case, 
 \begin{equation}\lambda_{ll}= \sum_{m=1}^{\infty}\frac{\Psi_l^2(x_m)}{\Psi_{m-1}^2(x_m)} \geq 1\end{equation}

 \begin{lemma}
  With the previous notations, there is a $\delta > \frac{5}{4}$ and a $c_0 = c_0(\delta)$ such that $f(c_0,\delta) \in (0,1)$ (see (38) for $f(c_0,\delta)$), with which
  \begin{equation}|\lambda_{kl}| \leq c\max\{k, l\}^{-\delta}, \ws\ws\ws \mbox{if} \ws \ws \ws k\neq l\end{equation}
 for a $c<c_0$ .\end{lemma}
 
 \proof
 
Let us suppose that $k < l$. We have to distinguish two cases: $\exists c>1$ such that $ck>l$, that is $k\sim l$, and $k<<l$. At first we will deal with the second case: with a $c > 1$ 
$$|\lambda_{kl}| \leq \sum_{m=1}^{ck}\frac{|\Psi_k(x_m)\Psi_l(x_m)|}{\Psi_{m-1}^2(x_m)} + \sum_{m=ck+1}^{\frac{1}{c}l} (\cdot) + \sum_{m=\frac{1}{c}l+1}^{cl}(\cdot) +\sum_{m=cl+1}^{\infty}(\cdot)$$ $$= S_1+S_2+S_3+S_4$$
In $S_2$ the first term of the nominator is exponentially small (see (21), (47)), that is by (20), (48), (52),(11)
$$S_2 \leq c \sum_{m=ck+1}^{\frac{1}{c}l}\frac{e^{-cQ\left(Q^{[-1]}(g(m))\right)}Q^{[-1]}(g(m))}{g^{\frac{1}{3}}(m)}\frac{1}{a_{\Psi_l}^{\frac{1}{4}}(a_{\Psi_l}-|x_m|)^{\frac{1}{4}}}$$
$$\leq \frac{c}{(Q^{[-1]}(g(l)))^{\frac{1}{2}}}\sum_{m=ck+1}^{\frac{1}{c}l}\frac{e^{-cQ\left(Q^{[-1]}(g(m))\right)}Q^{[-1]}(g(m))}{g^{\frac{1}{3}}(m)}$$
\begin{equation}\leq c\left(Q^{[-1]}(g(l))\right)^{-\frac{1}{2}} \leq c l^{-\frac{5}{4}},\end{equation}
where $c < \frac{c_0}{4},$ if $n_0$ is enough large.

$|\Psi_k(x_m)|$ is exponentially small in $S_3$, and there are $\sim l$ terms in it, furtheremore we can estimate $\frac{|\Psi_l(x_m)|}{|\Psi_{m-1}(x_m)|}$ by 1, so
\begin{equation}S_2 \leq c l \max_{m=\frac{1}{c}l+1}^{cl}|\Psi_k(x_m)|\frac{\left(Q^{[-1]}(g(cl))\right)^{\frac{1}{2}}}{g^{\frac{1}{6}}(cl)}\leq c le^{-g(cl)}\frac{\left(Q^{[-1]}(g(cl))\right)^{\frac{1}{2}}}{g^{\frac{1}{6}}(cl)}\leq c l^{-\frac{5}{4}},\end{equation}
where $c < \frac{c_0}{4},$ if $n_0$ is enough large.

In $S_4$, both the terms in the nominator  are exponentially small, that is
\begin{equation}S_4\leq c  \sum_{m=cl+1}^{\infty}\frac{e^{-cg(m)}Q^{[-1]}(g(m))}{g^{\frac{1}{3}}(m)} \leq c  l^{-\frac{5}{4}},\end{equation}
where $c < \frac{c_0}{4},$ if $n_0$ is enough large.

In $S_1$ we have to separate the "maximal" term:
$$S_1 \leq \sum_{m=1\atop m \neq k+1}^{ck}(\cdot) + \left|\frac{\Psi_l(x_{k+1})}{\Psi_k(x_{k+1})}\right| \leq \frac{c}{\left(Q^{[-1]}(g(k))\right)^{\frac{1}{4}}\left(Q^{[-1]}(g(l))\right)^{\frac{1}{4}}}$$ 
$$\times \sum_{m=1\atop m \neq k+1}^{ck}\frac{Q^{[-1]}(g(m))}{g^{\frac{1}{3}}(m)(a_{\Psi_k} - |x_m|)^{\frac{1}{4}}(a_{\Psi_l} - |x_m|)^{\frac{1}{4}}} $$ $$+ \frac{c}{\left(Q^{[-1]}(g(l))\right)^{\frac{1}{4}}}\frac{1}{(a_{\Psi_l} - |x_{k+1}|)^{\frac{1}{4}}}\frac{\left(Q^{[-1]}(g(k))\right)^{\frac{1}{2}}}{g^{\frac{1}{6}}(k)}$$
Because we deal with the $k<<l$ case we can estimate 
\begin{equation}|a_{\Psi_l} - |x_{k+1}|| > c a_{\Psi_l},\end{equation}
and so the second term, $S_{12}$, can be estimated as
$$S_{12} \leq c \frac{1}{\left(Q^{[-1]}(g(l))\right)^{\frac{1}{2}}}\frac{\left(Q^{[-1]}(g(k))\right)^{\frac{1}{2}}}{g^{\frac{1}{6}}(k)}.$$
As in Lemma 3, according to the behavior of the norm of orthogonal polynomials, we have to distinguish two cases: if the second factor is bounded in $k$, then by (11)
\begin{equation}S_{12} \leq c \frac{1}{\left(Q^{[-1]}(g(l))\right)^{\frac{1}{2}}}\leq c  l^{-\frac{5}{4}},\end{equation}
where $c < \frac{c_0}{8},$ if $n_0$ is enough large.

If the second factor is increasing, then also by (11)
 \begin{equation}S_{12} \leq c \frac{1}{g^{\frac{1}{6}}(l)}\leq c  l^{-\frac{5}{4}},\end{equation}
 where $c < \frac{c_0}{8},$ if $n_0$ is enough large.
 
 Now we have to deal with the first term of $S_1: S_{11}.$ 
 $$S_{11} \leq \frac{c}{\left(Q^{[-1]}(g(k))\right)^{\frac{1}{4}}\left(Q^{[-1]}(g(l))\right)^{\frac{1}{4}}}$$
 $$\times\sum_{m=1\atop m \neq k+1}^{ck}\frac{Q^{[-1]}(g(m))}{g^{\frac{1}{3}}(m)}\frac{1}{(a_{\Psi_k} - |x_m|)^{\frac{1}{4}}(a_{\Psi_l} - |x_m|)^{\frac{1}{4}}}$$
 As in (78)
 $$S_{11} \leq \frac{c}{\left(Q^{[-1]}(g(l))\right)^{\frac{1}{2}}}\sum_{m=1\atop m \neq k+1}^{ck}\frac{\left(Q^{[-1]}(g(m))\right)^{\frac{1}{2}}}{g^{\frac{1}{3}}(m)}\frac{\left(a_{\Psi_m}\right)^{\frac{1}{4}}}{\left(a_{\Psi_k}\right)^{\frac{1}{4}}}\frac{|x_m|^{\frac{1}{4}}}{(|x_{k+1}| - |x_m|)^{\frac{1}{4}}}$$
 \begin{equation}\leq \frac{c}{\left(Q^{[-1]}(g(l))\right)^{\frac{1}{2}}}\sum_{m=1\atop m \neq k+1}^{ck}m^{-\frac{5}{2}}m^{\frac{1}{4}}\leq \frac{c}{\left(Q^{[-1]}(g(l))\right)^{\frac{1}{2}}}\leq c  l^{-\frac{5}{4}},\end{equation}
 where $c < \frac{c_0}{8},$ if $n_0$ is enough large. Here we used (11), and the polynomially growing property of $g$ and $Q$.
 
 If $k \sim l$, then the second term is missing. In this case
 $$S_1 \leq \sum_{m=1\atop m \neq k+1,l+1}^{ck}(\cdot) + \left(\left|\frac{\Psi_l(x_{k+1})}{\Psi_k(x_{k+1})}\right| + \left|\frac{\Psi_k(x_{l+1})}{\Psi_l(x_{l+1})}\right|\right)= S_{11}+ S_{12},$$
 and as in (60), we can estimate by
 \begin{equation}|a_{\Psi_l} - |x_{k+1}|| > \frac{Q^{[-1]}(g(l))}{l}  ,\end{equation}
 and so by (11)
\begin{equation}S_{12} \leq c \frac{l^{\frac{1}{4}}}{g^{\frac{1}{6}}(l)}\leq c  l^{-\frac{5}{4}},\end{equation}
 where $c < \frac{c_0}{6},$ if $n_0$ is enough large.
 
 As in the previous calculation
 $$S_{11} \leq \frac{c}{\left(Q^{[-1]}(g(l))\right)^{\frac{1}{2}}}\sum_{m=1\atop m \neq k+1,l+1}^{ck}\frac{\left(Q^{[-1]}(g(m))\right)^{\frac{1}{2}}}{g^{\frac{1}{3}}(m)}\frac{|x_m|^{\frac{1}{4}}}{(|x_{k+1}| - |x_m|)^{\frac{1}{4}}}\frac{|x_m|^{\frac{1}{4}}}{(|x_{l+1}| - |x_m|)^{\frac{1}{4}}}$$
 \begin{equation}\leq \frac{c}{\left(Q^{[-1]}(g(l))\right)^{\frac{1}{2}}}\sum_{m=1\atop m \neq k+1,l+1}^{ck}m^{-\frac{5}{2}}m^{\frac{1}{4}}m^{\frac{1}{4}}\leq c  l^{-\frac{5}{4}},\end{equation}
 where $c < \frac{c_0}{6},$ if $n_0$ is enough large.
 
 When $k \sim l$, the estimations on $S_3$ and on $S_4$ are the same as in the previous case, and so the lemma is proved.
 
 Finally applying Lemma 2 with it's Corollary to the operator $\hat{A}^T\hat{A}$, and Lemma 4 to $\hat{A}^T$, we get that {\it Ran}$\hat{A}^T = l_2$, and so {\it Ker}$\hat{A}^ = \{0\}$, which proves the unicity of the solution.
 
 \subsection{Finite Section Method}
 
 As a consequence of invertibility we can apply the so-called "finite section method", which is a very natural (numerical) way to get the solution of the infinite equation $\hat{A}x = \hat{c}$. The process is the following: considering the system $Ax=b$, where $A$ is invertible, but not necessarily hermitian, we set 
 \begin{equation}A_{rn} = P_r AP_n, \ws \ws \ws \mbox{and}\ws\ws\ws b_{rn} = A^{*}_{rn}b, \end{equation}
 where $P_r$ and $P_n$ are projections, that is we can take $A_{rn}$ as it consists of the intersection of the first $r$ rows and the first $n$ columns of $A$, and $b_{rn}$ as the image of the cut vector $b_r$.
 Now we have to try to solve the equation
 \begin{equation}A^{*}_{rn}A_{rn}x_{rn}=b_{rn}.\end{equation}
 The convergence of this method is proved by K. Gr\"ochenig, Z. Rzesztonik and T. Strohmer \cite{gro}. To state the above-mentioned convergence theorem we need some notations. At first we have to note that the original paper works with the index class $\bf{Z}^d$, but without any modification we can apply the definitions and results to the index set $\bf{N}$.
 \begin{defi}
 We say that a matrix $A$ belongs to the Jaffard class ${\cal A}_s$, if it's elements $a_{kl}$, $k,l \in \bf{N}$ fulfil the following inequality:
 \begin{equation}|a_{k,l}| \leq C (1+ |k-l|)^{-s} \ws\ws\ws \forall k,l \in \bf{N},\end{equation}
 where $C$ is an absolute constant. The norm in the Jaffard class is $\|A\|_{{\cal A}_s} = \sup_{k,l \in \bf{N}}|a_{k,l}| (1+ |k-l|)^{s}$. \end{defi}
 
 \note
 
 Let us denote by $\sigma(A^{*}A)$ the spectrum of $A^{*}A$, and by $\lambda_{-} = \min \sigma(A^{*}A)$.
 
 So the (simplified version of the) theorem is the following \cite{gro}, Th. 16:
 
 \noindent {\it {\bf Theorem C}
 
 Let $A \in {\cal A}_s$ with an $s>1$, and $Ax=b$ is given, where $b \in l^2$, and $A$ is invertible on $l^2$. Consider the finite sections
 $$A^{*}_{rn}A_{rn}x_{rn}=b_{rn}.$$
 Then, for every $n$ there exists an $R(n)$ (depending on $\lambda_{-}$ and $s$), such that $x_{r(n)n}$ converges to $x$ in the norm of $l^2$, for every choice $r(n) \geq R(n).$} 
 
 Because $\hat{c}$ is in $l^2$, and $\hat{A}$ is invertible on $l^2$, for the convergence of the finite section method we have to prove, that $\hat{A} \in {\cal A}_s$ with an $s>1$.
 
 \begin{lemma} Let $\hat{A}$ be as in (30), and let us suppose the assumption of Definition 4, then there exists an $s>1$ such that $\hat{A} \in {\cal A}_s$ \end{lemma}
 
 \proof
 
 Because in our matrix $\hat{A}$ the dominant elements are under the principal diagonal, ($|a_{k,k-1}|=1$), we have to shift the indices, that is we have to prove that there is an absolute constant $C$, and there is an $s>1$ such that 
 \begin{equation}|a_{k,l-1}| = \left|\frac{\Psi_{l-1}(x_k)}{\Psi_{k-1}(x_k)}\right| \leq C (1+ |k-l|)^{-s}.\end{equation}
 Using (49) and (52), we get that
 $$|a_{k,l-1}|\leq c \frac{\left(Q^{[-1]}(g(k))\right)^{\frac{1}{2}}}{\left(g(k)\right)^{\frac{1}{6}}}\frac{1}{\left(Q^{[-1]}(g(l))\right)^{\frac{1}{4}}}\frac{1}{\left|a_l - |x_k|\right|^{\frac{1}{4}}},$$
 where $c$ is an absolute constant. 
 
 Considering that $x_k \sim a_{\Psi_k}$, we have to distinguish some cases: 
 
 \noindent a): if $l \sim k$, but $|a_{\Psi_l} - |x_k||$ is not too small, or if $k << l$, then with some $c_1 \neq c_2$
 \begin{equation}D_{lk} = \left|a_l - |x_k|\right|^{\frac{1}{4}} \geq c\left|c_1Q^{[-1]}(g(l) - c_2 Q^{[-1]}(g(k)\right|^{\frac{1}{4}}\geq c \left(Q^{[-1]}(g(l)\right)^{\frac{1}{4}}.\end{equation}
 
 \noindent b): if $l<<k$, then
 \begin{equation}D_{lk} \geq c \left(Q^{[-1]}(g(k)\right)^{\frac{1}{4}}.\end{equation}
 
 \noindent c): if $l \sim k$, and $x_k$ is close to $a_l$, recalling the estimation on the distance of two consecutive maximum points of $\Psi_k$-s, (see (61))
 \begin{equation}D_{lk} \geq c \left(\frac{Q^{[-1]}(g(l))}{l}\right)^{\frac{1}{4}},\end{equation}
 where we used the polinomially growing property of $Q$ and $g$.
 
 So in case a), when $l\sim k$ we get that
 \begin{equation}|a_{k,l-1}| \leq c \frac{1}{\left(g(l)\right)^{\frac{1}{6}}} \leq C (1+ |k-l|)^{-s},\end{equation}
 for all $s>1$.
 
 Also in case a), when $k<<l$
 \begin{equation}|a_{k,l-1}| \leq c \left(\frac{Q^{[-1]}(g(k))}{Q^{[-1]}(g(l))}\right)^{\frac{1}{2}}\frac{1}{\left(g(k)\right)^{\frac{1}{6}}}\end{equation}
 Here, as in Lemma 3 ((53),(54)), we have to distinguish two cases: in the first case, according to (11) 
 \begin{equation}|a_{k,l-1}| \leq c \frac{1}{\left(Q^{[-1]}(g(l)\right)^{\frac{1}{2}}} \leq c l^{-\frac{5}{4}}\leq C (1+ |k-l|)^{-s},\end{equation}
 with $s = \frac{5}{4}$. In the second case, according to (11) again
 \begin{equation}|a_{k,l-1}| \leq c \frac{1}{\left(g(l)\right)^{\frac{1}{6}}} \leq C (1+ |k-l|)^{-s},\end{equation}
 and $s>\frac{5}{4}$.
 
 In case b), also by (11)
 \begin{equation}|a_{k,l-1}| \leq c \left(\frac{Q^{[-1]}(g(k))}{Q^{[-1]}(g(l))}\right)^{\frac{1}{4}}\frac{1}{\left(g(k)\right)^{\frac{1}{6}}}\leq c \frac{\left(Q^{[-1]}(g(k))\right)^{\frac{1}{4}}}{\left(g(k)\right)^{\frac{1}{6}}}\leq C (1+ |k-l|)^{-s}\end{equation}
 with $s = \frac{5}{4}$.
 
 In case c), by (10),(11)
 \begin{equation}|a_{k,l-1}| \leq c \frac{l^{\frac{1}{4}}}{\left(g(l)\right)^{\frac{1}{6}}}\leq c\frac{\left(Q^{[-1]}(g(l)\right)^{\frac{1}{4}}}{\left(g(l)\right)^{\frac{1}{6}}} \leq C (1+ |k-l|)^{-s},\end{equation}
 for all $s>1$, wich proves the lemma.

 \subsection{Convergence}
 
 As it turned out in the introduction, the required form of the elements of the dual space is the following:
 $$\varphi_m^{*}= \frac{\varphi_m-\sum_{k=1}^{\infty}a_{km}\Psi_{k}}{v^2},$$
 which implies that we have to deal with the convergence of the series in the nominator, and we have to give some estimations on the order of the zeros of the nominator.
 
 \begin{lemma}Let $\{a_{km}\}_{k=1}^{\infty}$ be an $l_2$-solution of (30), then 
 \begin{equation}\left|\sum_{k=1}^na_{km}\Psi_k(x)\right| \leq  c \frac{Q^{\frac{1}{6}}(x)}{(1+x^2)^{\frac{1}{4}}}\ws \ws n\in\mathbb{N},\end{equation}
 and the sum: $\sum_{k=1}^{\infty}a_{km}\Psi_k(x)$ is convergent in every $x \in\mathbb{R}$.
 \end{lemma}
 
 \proof
 
 At first we will show that the partial sum $\sum_{k=1}^na_{km}\Psi_k(x)$ can be estimated by a function wich grows at most polynomially on $\mathbb{R}$. Using Cauchy-Schwarz's inequality, we have to estimate
 $$ s_n(x) = \sum_{k=1}^n\Psi_k^2(x)=\Psi_{j(x)}^2(x)+ \sum_{k=1}^{\frac{1}{c}j(x)}\Psi_k^2(x) + \sum_{\frac{1}{c}j(x) < k < cj(x) \atop k \neq j(x)}(\cdot) + \sum_{k=cj(x)}^n(\cdot)$$
 \begin{equation}=\Psi_{j(x)}^2(x)+ S_1 + S_2 + S_3 ,\end{equation}
 where $c>1$, and $j(x)$ means that index, for which the maximum point of $|\Psi_{j(x)}(x)|$ is the closest to $x$.  (Because $\Psi_k^2$-s are even, we can work on the positive part of the real line.) Hence, because the $n^{th}$ orthonormal polynomial $p_n$ attains its maximum around $a_n$, according to (20) and (49) we obtain that
 \begin{equation}\Psi_{j(x)}^2(x) \leq \|p_{g(j(x))}\|_{\infty}^2 \sim \left(g(j(x))\right)^{\frac{1}{3}}a_{g(j(x))}^{-1} \sim \frac{Q^{\frac{1}{3}}(x)}{x}\end{equation}
 By (48)
 \begin{equation} S_2 \leq c \sum_{\frac{1}{c}j(x) < k < cj(x) \atop k \neq j(x)}\frac{1}{\sqrt{a_{\Psi_k}}\sqrt{|a_{\Psi_k}-|x||}}.\end{equation}
 Taking into consideration the properties of $g$ and $Q$, we can estimate the difference under the square root as 
 $$|a_{\Psi_k}-|x||\geq c |a_{\Psi_k}-a_{\Psi_{j(x)}}| \geq c \left(Q^{[-1]}(g(\cdot))\right)^{'}(j(x))|k-j(x)|$$
 \begin{equation}\geq c \frac{Q^{[-1]}(g(j(x)))}{j(x)}|k-j(x)|,\end{equation}
 so 
 $$ S_2 \leq c \frac{\sqrt{j(x)}}{Q^{[-1]}(g(j(x))}\sum_{\frac{1}{c}j(x) < k < cj(x) \atop k \neq j(x)}\frac{1}{\sqrt{|k-j(x)|}}$$
 \begin{equation}\leq c \frac{j(x)}{Q^{[-1]}(g(j(x))} \leq c \frac{g^{[-1]}(Q(x))}{x}.\end{equation}
 As in Lemma 3 we collected the exponentially small terms in $S_1$, thus 
 \begin{equation} S_1 \leq c j(x)e^{-\frac{Q(x)}{2}} \leq c g^{[-1]}(Q(x))e^{-\frac{Q(x)}{2}}.\end{equation}
 In $S_3$, $x \sim a_{\Psi_{j(x)}}$, is far away from $a_{\Psi_k}$, so $|a_{\Psi_k}-a_{\Psi_{j(x)}}|  \geq c a_{\Psi_k}$, that is 
 \begin{equation} S_3 \leq c \sum_{k=cj(x)}^n \frac{1}{a_{\Psi_k}} \leq c\sum_{k=cj(x)}^n \frac{1}{Q^{[-1]}(g(k))}\end{equation}
 We can estimate this sum by
 \begin{equation}\int_{cj(x)}^{\infty}\frac{1}{Q^{[-1]}(g(y))}dy \leq \int_{cQ^{[-1]}(g(j(x))}^{\infty}\frac{1}{z}\frac{g^{[-1]}(Q(z))}{z}dz , \end{equation}
where we used the properties of $Q$ and $g$ again. By (10) we have that
$$S_3 \leq \sup_{z\geq cQ^{[-1]}(g(j(x))}\frac{g^{[-1]}(Q(z))}{z^{1-\varepsilon}}\int_{cQ^{[-1]}(g(j(x))}^{\infty}\frac{1}{z^{1+\varepsilon}}dz  \leq c \frac{j(x)}{Q^{[-1]}(g(j(x))}$$
\begin{equation}\leq c \frac{g^{[-1]}(Q(x))}{x}.\end{equation}
In the estimations above, we can replace $x$ in the denominator by $\sqrt{1+x^2}$ and so we haven't problem at zero. Collecting our estimations, if $g^{[-1]}(Q(x))$ is less then $Q^{\frac{1}{3}}(x)$ (see (11)), we obtain that
\begin{equation}\sum_{k=1}^n\left|a_{km}\Psi_k(x)\right|\leq c \|\{a_{km}\}\|_2 \frac{Q^{\frac{1}{6}}(x)}{(1+x^2)^{\frac{1}{4}}},\end{equation}
which gives uniform convergence, if $B < 3$, and the second statement of the lemma otherwise. 

\remark

The same computation yields that
\begin{equation}\left|\sum_{k=n}^{\infty}a_{km}\Psi_k(x)\right|\leq c \sqrt{\sum_{k=n}^{\infty}a_{km}^2} \frac{Q^{\frac{1}{6}}(x)}{(1+x^2)^{\frac{1}{4}}},\end{equation}
which means that $\sum_{k=1}^{\infty}a_{km}\Psi_k(x)$ tends to a function $f(x)$ locally uniformly on $\mathbb{R}$.

Let $\varphi_l$ be an element of the system (16). Considering that $\varphi_l$ is a weighted polynomial with an exponential weight, we can immediately get the following

\corollary

There exists a function $g \in L^1({\bf R})$ such that
\begin{equation}\left|\sum_{k=1}^na_{km}\Psi_k(x)\varphi_l(x)\right| \leq  g(x)\ws \ws n\in\mathbb{N}.\end{equation}

To state the following lemma we need some notations. Let $S_j := (x_j- \delta(x_j),x_j+ \delta(x_j)$ is a ball around $x_j$ such that $x_i \notin S_j$ if $i\neq j$. And let
$$\sigma_n(x) = \sigma_{n,m}(x)= \sum_{k=0}^n\left(1-\frac{l_k}{n+1}\right)a_{km}\Psi_k(x)$$
be the $n^{th}$ Cesaro mean of the Fourier series with respect to $\{p_n(w)\}_{n=0}^{\infty}$ of $S = \sum_{k=1}^{\infty}a_{km}\Psi_k(x)$, where $\{a_{km}\} \in l_2$ is the solution of (30). With these notations we have

\begin{lemma} Supposing (10)
\begin{equation}|\sigma_n(x) -\sigma_{\nu}(x)| = O\left(\frac{1}{n^{\gamma}}\right) \ws \ws \mbox{if}\ws\ws x\in S_j,\ws\ws n(j)< n <\nu, \ws\ws j=1,2 \dots \end{equation}
\end{lemma}

\proof

Let $x$ be in $S_j$.
$$|\sigma_n(x) -\sigma_{\nu} (x)| \leq \left|\sum_{0 \leq l_k \leq n}\left(\frac{1}{n+1}-\frac{1}{\nu+1}\right)l_ka_{km}\Psi_k(x)\right|$$ \begin{equation} + \left|\sum_{n < l_k \leq \nu}\left(1 - \frac{l_k}{\nu+1}\right)a_{km}\Psi_k(x)\right|= (*)\end{equation}

Let us denote by $k(x)$ that index for which the maximum point of $\Psi_{k(x)}$ is the nearest to $x$. If $x$ is around $a_{l_{k(x)}}$, then $x \sim Q^{[-1]}(l_{k(x)})$, that is $Q(x) \sim l_{k(x)}$. So if $n>N=N(j)$ (eg $cQ(x_j)<n $), then $cl_{k(x)} <n$. Let us assume now that $n$ is enough large:
$$(*) \leq \left|\sum_{0 \leq l_k \leq n \atop l_k \neq l_{k(x)}}\left(\frac{1}{n+1}-\frac{1}{\nu+1}\right)l_ka_{km}\Psi_k(x)\right| $$ $$+ \left|\left(\frac{1}{n+1}-\frac{1}{\nu+1}\right)l_{k(x)}a_{k(x)m}\Psi_{k(x)}(x)\right|$$ \begin{equation}+ \left|\sum_{n < l_k \leq \nu}\left(1 - \frac{l_k}{\nu+1}\right)a_{km}\Psi_k(x)\right|= S_1 + M + S_2\end{equation}

Let us recall (52) and at first we will deal with $M$.
\begin{equation}M \leq \frac{|a_{k(x)m}|}{n}l_{k(x)}\left(l_{k(x)}\right)^{\frac{1}{6}}\left(a_{l_{k(x)}}\right)^{-\frac{1}{2}} \leq c\frac{Q^{\frac{7}{6}}(x_j)}{\sqrt{x_j}}\frac{1}{n} = O\left(\frac{1}{n}\right)\end{equation}

We can handle $S_2$ as $S_3$ in Lemma 3, that is 
$$S_2 \leq \|\{a_{km}\}_{k=n}^{\infty}\|_2\left(\sum_{n < l_k \leq \nu}\left(1 - \frac{l_k}{\nu+1}\right)^2\Psi_k^2(x)\right)^{\frac{1}{2}}$$
$$\leq c\|\{a_{km}\}_{k=n}^{\infty}\|_2\left(\sum_{n < l_k \leq \nu}\frac{1}{a_{\Psi_k}}\right)^{\frac{1}{2}}$$ \begin{equation} \leq c \|\{a_{km}\}_{k=n}^{\infty}\|_2\left(\int_{cg^{[-1]}(n)}^{\infty}\frac{1}{Q^{[-1]}(g(y))}dy\right)^{\frac{1}{2}}= o\left(\sqrt{\frac{g^{[-1]}(n)}{Q^{[-1]}(n)}}\right)\end{equation}
We have to decompose $S_1$ to three parts. Using the Cauchy-Schwarz inequality again we obtain that
$$S_1 \leq \frac{c}{n}\left(\sum_{0 \leq l_k \leq \frac{1}{c}l_{k(x)}}l_k^2\Psi_k^2\right)^{\frac{1}{2}}+ \frac{c}{n}\left(\sum_{\frac{1}{c}l_{k(x)}< l_k \leq c l_{k(x)}}l_k^2\Psi_k^2\right)^{\frac{1}{2}}$$ \begin{equation}+\frac{c}{n}\left(\sum_{cl_{k(x)}< l_k \leq n}l_k^2\Psi_k^2\right)^{\frac{1}{2}}=S_{11}+S_{12}+S_{13}\end{equation}
Henceforward $S_{11}$ is the collection of the exponentially small terms, that is
\begin{equation} S_{11} \leq \frac{c}{n} \sqrt{Q(x_j)e^{-cQ(x_j)}} = O\left(\frac{1}{n}\right)\end{equation}
Applying also the same chain of ideas as in Lemma 3, we obtain that

$$ S_{12}\leq \frac{c}{n}l_{k(x)}\left(\sum_{\frac{1}{c}k(x)< k \leq c k(x)}\Psi_k^2\right)^{\frac{1}{2}} $$ 

\begin{equation} \leq \frac{cl_{k(x)}}{n}\sqrt{\frac{g^{[-1]}(Q(x))}{x}}\leq c Q(x_j)\sqrt{\frac{g^{[-1]}(Q(x_j))}{x_j}}\frac{1}{n} = O\left(\frac{1}{n}\right) \end{equation}
Similarly
$$S_{13} \leq \frac{c}{n}\left(\sum_{k=ck(x)}^nl_k^2\Psi_k^2\right)^{\frac{1}{2}} \leq \frac{c}{n}\left(\sum_{k=ck(x)}^n\frac{l_k^2}{a_{\Psi_k}}\right)^{\frac{1}{2}} $$
$$\sum_{k=ck(x)}^n\frac{l_k^2}{a_{\Psi_k}} \leq c \int_{cg^{[-1]}(Q(x))}^{g^{[-1]}(n)}\frac{g^2(y)}{Q^{[-1]}(g(y))}dy
$$ $$\leq c \int_{cx}^{Q^{[-1]}(n)}\frac{g^{[-1]}(Q(z))}{z^2}Q^2(z)dz \leq c g^{[-1]}(n)\int_{cx}^{Q^{[-1]}(n)}\frac{Q^2(z)}{z^2}dz $$
Applying \cite{lelu} 5.4. we can estimate $\frac{Q(z)}{z}$ by $\frac{1}{A}Q^{'}(z) $, where $A>1$ is in the definition of Freud weights. So
$$\int_{cx}^{Q^{[-1]}(n)}\frac{Q^2(z)}{z^2}dz  \leq \frac{1}{2A}\int_{cx}^{Q^{[-1]}(n)}2Q^{'}(z)Q(z)\frac{1}{z}dz$$
With an integration by parts we get that
$$\int_{cx}^{Q^{[-1]}(n)}\frac{Q^2(z)}{z^2}dz \leq \frac{1}{2A-1}\left(\frac{n^2}{Q^{[-1]}(n)}- c\frac{Q^2(x)}{x}\right) \leq c \frac{n^2}{Q^{[-1]}(n)},$$ if $n$ is large enough. Summarizing the calculations of the previous lines we obtain that
 \begin{equation} S_{13} \leq c \sqrt{\frac{g^{[-1]}(n)}{Q^{[-1]}(n)}}\end{equation}
 Hence these estimations yield that if (10) fulfils, then $|\sigma_n(x) -\sigma_{\nu}(x)| = O\left(\frac{1}{n^{\gamma}}\right)$.
 
 This lemma showes the order of the roots of $\varphi_m- \sum_{j=1}^{\infty}a_{jm}\Psi_j$ at $x_j$-s, that is applying the classical theorem of S. N. Bernstein (see eg. \cite{nat}), and taking into consideration Lemma 7 and its corollary as well, we get the following
 
 \newpage
 
 \corollary
 
 Let $x\in S_j$, then
 \begin{equation}\left|\varphi_m(x)- \sum_{j=1}^{\infty}a_{jm}\Psi_j(x)\right| \leq |x-x_j|^{\gamma}h(x), \end{equation}
 where $h(x)$ is independent of $j$, is continuous and it grows polynomially with $x$.

  For  the final computations let us prove our last lemma, which  follows the same chain of ideas as Lemma 1.1 of J. Szabados \cite{joska}:
 
 \begin{lemma}
 
 Let $m_j, \varrho\geq 0, \varepsilon > 0, x_j$ be as in Definition 2, with properties (2),and let 
 $$\hat{v}(x) = \prod_{j=1}^{\infty}\left|1-\frac{x}{x_j}\right|^{m_j} \ws\ws\mbox{and} \ws\ws  \hat{v}_k(x) = \prod_{1\leq j<\infty \atop j \neq k}\left|1-\frac{x}{x_j}\right|^{m_j}$$
 Then
 \begin{equation}\hat{v}(x) \leq e^{c|x|^{\varrho + \varepsilon}}, \ws \ws x\in\mathbb{R}, \ws \ws \varepsilon > 0\end{equation}
 and
 \begin{equation}\hat{v}(x) \geq e^{-c|x|^{\varrho + \varepsilon}}, \ws \ws \mbox{for} \ws x \in \mathbb{R}\setminus \cup_{j=1}^{\infty}\left(x_{j}-\frac{m_{j}}{|x_{j}|^{\varrho + \varepsilon}},x_{j}+\frac{m_{j}}{|x_{j}|^{\varrho + \varepsilon}}\right),\end{equation}
 furtheremore
 \begin{equation}\hat{v}_k(x) \geq e^{-c|x|^{\varrho + \varepsilon}}, \ws \ws \mbox{for} \ws x \in \left(x_{k}-\frac{m_{k}}{|x_{k}|^{\varrho + \varepsilon}},x_{k}+\frac{m_{k}}{|x_{k}|^{\varrho + \varepsilon}}\right),\end{equation}
 where $c>0$ depends on $\hat{v}$ and $\varepsilon$, and if $a>b$, then $[a,b]= \emptyset$.\end{lemma}
 
\remark

If eg.  $x_j = j^{\nu}, \nu > 0$, then $\varrho=\frac{1}{\nu}$, and if $x_j = 2^{j}$, then $\varrho=0$.

The proof follows the steps of the proof of Lemma 1.1 in \cite{joska}.

\proof

Let us denote by
$$N(x)= \sum_{|x_k|<|x|}m_k $$
According to (4), $N(x) \leq c(\varepsilon)|x|^{\varrho+\varepsilon},$ for all $\varepsilon>0$.
$$\hat{v}(x) \leq \prod_{|x_k|<|x|}\left|1-\frac{x}{x_k}\right|^{m_k}\prod_{|x_k|\geq|x|\atop x_kx<0}\left(1+\left|\frac{x}{x_k}\right|\right)^{m_k}= \hat{v}_1(x)\hat{v}_2(x)$$
As in \cite{joska}, 
$$\hat{v}_1(x) \leq \prod_{|x_k|<|x|}\left(2\left|\frac{x}{x_k}\right|\right)^{m_k}\leq (2|x|)^{N(x)}\left(\frac{\sum_{|x_k|<|x|}\frac{m_k}{|x_k|^{\varrho+\varepsilon}}}{N(x)}\right)^{\frac{N(x)}{\varrho+\varepsilon}}$$
$$\leq \left(\frac{c|x|^{\varrho+\varepsilon}}{N(x)}\right)^{\frac{N(x)}{\varrho+\varepsilon}}\leq e^{c|x|^{\varrho+\varepsilon}}.$$
$$\hat{v}_2(x) \leq e^{\sum_{|x_k|\geq|x|\atop x_kx<0}m_k\log\left(1+\left|\frac{x}{x_k}\right|\right)}\leq e^{|x|^{\varrho+\varepsilon}\sum_{k=1}^{\infty}\frac{m_k}{|x_k|^{\varrho+\varepsilon}}} \leq e^{c|x|^{\varrho+\varepsilon}}.$$
For the lower estimation, as in \cite{joska}, we devide our product to three parts: if $x \neq x_j, j=1,2,\dots$
$$\hat{v}(x)=\prod_{|x_j|< |x|}\left|1-\frac{x}{x_j}\right|^{m_j}\prod_{|x|<|x_j|\leq 2|x|}(\cdot)\prod_{|x_j|> 2|x|}(\cdot)=P_1P_2P_3$$
As $P_1\geq \prod_{\frac{x}{x_j}>1}(\cdot), P_2\geq \prod_{1<\frac{x}{x_j}<2}(\cdot), P_3\prod_{0<\frac{x}{x_j}<\frac{1}{2}}(\cdot)$ the computations are the same as in \cite{joska}, so we omit the details.

Also the same computation implies (126).

Now we are in the position to prove the theorem.

 \subsection{Proof of the Theorem}

The properties of $g$ imply that $\varrho < 1$ in the definition of $v$, so it is obvious from (123) and the definition of Freud weight, that there exists a $\mu$, such that with arbitrary $d>0$ there is a $v:= v_{X,M,\mu,d}$ with which $\varphi_k v \in L^p$. According to (125) and (126) to $c = c(\mu)$ we can choose a $d>0$ such that $v_{X,M,\mu,d} > ce^{k|x|^{\varrho+\mu}}$ with some $k>0$ on $\mathbb{R}\setminus \cup_{j=1}^{\infty}\left(x_{j}-\frac{m_{j}}{|x_{j}|^{\varrho + \varepsilon}},x_{j}+\frac{m_{j}}{|x_{j}|^{\varrho + \varepsilon}}\right)$, and the same fulfils on $v_k= \hat{v}_k(x)e^{d|x|^{\varrho+\mu}}$ on the interval $\left(x_{k}-\frac{m_{k}}{|x_{k}|^{\varrho + \varepsilon}},x_{k}+\frac{m_{k}}{|x_{k}|^{\varrho + \varepsilon}}\right)$.

 Let
 \begin{equation}\varphi_m^{*} = \frac{1}{v^2}\left(\varphi_m- \sum_{j=1}^{\infty}a_{jm}\Psi_j\right)\ws \ws m=1,2,\dots,\end{equation}
 where $\{a_{jm}\}$ is a solution of (30).
 We will show that $\{\varphi_m^{*}\}_{m=1}^{\infty}$ is a system in $L^q_v$ which is biorthonormal with respect to $\{\varphi_m\}_{m=1}^{\infty}\subset L^p_v$.
 
  According to Lemma 7, the series in (101) is convergent in some sense, that is the definition of $\varphi_m^{*}$ is clear, and applying the Corollary after Lemma 7, by Lebesgue's theorem we can integrate term by term in the followings:
 $$\int_{\mathbb{R}}\varphi_m^{*}\varphi_k v^2= \int_{\mathbb{R}}\frac{1}{v^2}\left(\varphi_m- \sum_{j=1}^{\infty}a_jm\Psi_j\right)\varphi_k v^2$$ $$=\int_{\mathbb{R}}\varphi_m\varphi_k   - \sum_{j=1}^{\infty}a_{jm}\int_{\mathbb{R}}\Psi_j\varphi_k  = \delta_{m,k} + 0,$$
 where we used the orthonormality of the original system, which was the weighted othonormal polynomials.
 
 So the only thing we have to prove that $\varphi_m^{*}$ is in $L^q_v$. Let $S^{*}_j=S_j\cap (x_{j}-\frac{m_{j}}{x_{j}^{\varrho + \varepsilon}},x_{j}+\frac{m_{j}}{x_{j}^{\varrho + \varepsilon}})$, and $\delta > 0$ (see 5)) is fixed.  Thus
 $$\int_{\mathbb{R}}\left|\frac{\varphi_m- \sum_{j=1}^{\infty}a_{jm}\Psi_j}{v}\right|^q= \sum_{j=1}^{\infty}\int_{x \in S^{*}_j}(\cdot)+ \int_{x \in (\mathbb{R}\setminus \ \cup_{j=1}^{\infty}S^{*}_j)}(\cdot)$$
 $$\leq c \sum_{j=1}^{\infty}\int_{x \in S^{*}_j}\left|h(x)x^{m_j}e^{-k|x|^{\varrho+\mu}} \right|^q\left||x-x_j|^{\gamma-m_j}\right|^q +$$ $$ \int_{x \in (\mathbb{R}\setminus \cup_{j=1}^{\infty}S^{*}_j)}\left|k(x)e^{-k|x|^{\varrho + \mu}}\right|^q =(*),$$
 where $h(x)$ is as in (122), and according to Lemma 7, $k(x)$ grows polynomially. Hence by Lemma 9, we can estimate $(*)$ on the whole real line with an integral of a function grows polynomially, times an exponentially small factor, that is
 $$\|\varphi_m^{*}\|_q \leq \left(\int_\mathbb{R}\left|k_1(x)e^{-c|x|^{\varrho + \mu}}\right|^q\right)^{\frac{1}{q}},$$
 where $k_1(x)$ depends only on $m$ and $Q(x)$, so the $q-$norm of $\varphi_m^{*}$ is bounded if $m_j-\gamma < \frac{1}{q}$, so the dual system in $L^q_v$, when $p$ fulfils the inequalities in the theorem.
 
 For completeness we have to prove, that if for a $g \in L^q_{v_X,M}$ (where $\frac{1}{p}+\frac{1}{q}=1$) $g(\varphi_k)= \int_{\mathbb{R}}g\varphi_kv^2 = 0 , k\in \mathbb{N},$ then $g=0.$ The comleteness of the original system implies that $g$ has to be in form: 
 $$g= \frac{1}{v^2}\sum_{j=1}^{\infty}b_j\Psi_j,$$
 and as $g \in L^q_{v_{X,M,d,\mu}}$, $\int_{\mathbb{R}}|gv|^q$ must be finite. By the properties of $v$, and recalling that $ \Psi_j = p_{l_j}w$, the integral on $\mathbb{R}\setminus \cup_jS_j^{*}$ is finite, so we have to deal with the integral around the roots of $v$, that is $\sum_{j=1}^{\infty}\int_{S_j^{*}}\left|\frac{1}{v}\sum_{j=1}^{\infty}b_j\Psi_j\right|^q$ has to be finite. Together with the assumption: $p< \inf_{m_j <1}\frac{1}{1-m_j}$, it means that
 \begin{equation}\left(\sum_{j=1}^{\infty}b_j\Psi_j\right)(x_k) = 0 \hspace{1cm} k=1,2,\dots \end{equation}
 So as in (15), we got a homogene linear equation system:
 \begin{equation} A b = 0,\end{equation}
 where $A$ is the same infinite matrix as in (15). Introducing $\hat{A}$, etc, according to 3.1.2, the homogene equation has the only solution in $l_2$: $b_j = 0, j=1,2 \dots$, that is $g=0$.
 
 \noindent {\bf Final Remarks}
 
 (A) If somebody doesn't take care on the range of the operator $A$, then, because in our case on the right hand side of the equation there is a fast convergent vector,  to get some solution of the equation $Aa_m=c_m$, it is enough to apply Toeplitz's theorem, and so it is not necessary to guarantee a not too small element in every rows. That is the proof of Lemma 1 ensures a good omission system for arbitrary point systems.
 
 (B) The aim of this paper was to show the existence of a "good" point- and a "good" omission system with some assumptions on the functions $Q$ and $g$. We chose a rather comfortable one. More precisely, our calculations show that besides (10), which needs for convergence, it is enough to assume
 
 for solvability:
 $$g(x) > x^{\mu}, \ws\ws\ws \mu>\frac{15}{2},$$
 $\frac{g^{[-1]}(x)}{(Q^{[-1]}(x))^{1-\varepsilon}}$ is strictly decreasing for a $\varepsilon>0$,
 $$x^{\delta}\max\left\{\frac{x^{\frac{1}{4}}}{g^{\frac{1}{6}}(x)};\frac{1}{(Q^{[-1]}(x))^{\frac{1}{2}}}\right\} \to 0, \ws\ws\ws\mbox{where}\ws\ws\ws \delta > \frac{5}{4}$$
 
 for unicity:
 
 $$x^{\delta}\max\left\{\frac{1}{g^{\frac{1}{6}}(x)};\frac{1}{(Q^{[-1]}(x))^{\frac{1}{2}}}\right\} \to 0, \ws\ws\ws\mbox{where}\ws\ws\ws \delta > \frac{5}{4},$$
 and
 $$x^{\nu}\frac{(Q^{[-1]}(g(x)))^{\frac{1}{4}}}{g^{\frac{1}{6}}(x)}\to 0, \ws\ws\ws\mbox{where}\ws\ws\ws \nu > \frac{3}{4},$$
 
 for the convergence of finite section method:
 $$x^{\kappa}\max\left\{\frac{(Q^{[-1]}(g(x)))^{\frac{1}{4}}}{g^{\frac{1}{6}}(x)};\frac{1}{(Q^{[-1]}(x))^{\frac{1}{2}}}\right\} \to 0, \ws\ws\ws\mbox{where}\ws\ws\ws \kappa > 1$$

\end{document}